\documentclass{article}

\usepackage{amssymb} 
\usepackage{amsmath} 
\usepackage{theorem} 

\usepackage{epsfig}

\theorembodyfont{\itshape} 

\newtheorem{theorem}{Theorem}[section]
\newtheorem{lemma}[theorem]{Lemma}
\newtheorem{corollary}[theorem]{Corollary}
\newtheorem{proposition}[theorem]{Proposition}

\newenvironment{proof}{{\noindent \bfseries\itshape Proof: }}{\hfill$\Box$} 

\newenvironment{proofofproposition}{{\noindent \bfseries\itshape Proof of Proposition }}{\hfill$\Box$} 
\newenvironment{proofoftheorem}{{\noindent \bfseries\itshape Proof of Theorem }}{\hfill$\Box$}

\def\e{\varepsilon}
\def\p{\prime}
\def\Si{\varSigma}
\def\Ga{\varGamma}
\def\vect#1{\mbox{\boldmath $#1$}} 
\def\nekofrac#1#2{\mbox{\footnotesize$\displaystyle\frac{#1}{#2}$}}
\def\wfrac#1#2{\frac{\,#1\,}{\,#2\,}} 
\DeclareFontFamily{U}{rsfs}{\skewchar\font"7F}
\DeclareFontShape{U}{rsfs}{m}{n}{
	<-6> rsfs5
	<6-8> rsfs7
	<8-> rsfs10
	}{}
\DeclareMathAlphabet{\mathscr}{U}{rsfs}{m}{n}
\usepackage{color}

\begin{document}

\title{\bf Conformal invariance of the writhe of a knot}

\author{R. ~Langevin \\{\small Institut de Math\'ematiques de Bourgogne, Universit\'e de Bourgogne} \and J. ~O'Hara\footnote{This joint work was done during the second author's stay in Dijon which was supported by JSPS.} \\{\small Department of Mathematics, Tokyo Metropolitan University}}

\maketitle

\begin{abstract}
We give a new proof of the conformal invariance of the writhe of a knot from a conformal geometric vewpoint. 
\end{abstract}

\medskip
{\small {\it Key words and phrases}. writhe, conformal invariant.}

{\small 2000 {\it Mathematics Subject Classification.} Primary 57M25, Secondary 53A30.}



\section{Introduction} 
Suppose $K$ is a framed knot, i.e. there is a unit normal vector field $\vect e_2$ along $K$. 
A $2$-component link $K\cup K+\e\vect e_2$ $(|\e|\ll1)$ can be considered a closed ribbon. 
Let $\textsl{Lk}$ and $\textsl{Tw}$ be the linking number and the total twist of $K\cup K+\e\vect e_2$, and $\textsl{Wr}$ the writhe of $K$. 
Then we have 
\begin{equation}\label{Lk=Wr+Tw}
\textsl{Lk}=\textsl{Wr}+\textsl{Tw}. 
\end{equation}
When the knot $K$ has nowhere vanishing curvature and $\vect e_2$ is the principal normal vector, $\textsl{Lk}$ is called the {\sl self-linking number} of the knot and denoted by $\textsl{Sl}$. 
The equation (\ref{Lk=Wr+Tw}) was proved in this case in \cite{Ca1, Ca2, Ca3, Po1}, and in \cite{Wh} in general. 
It plays an important role in the application of the knot theory to molecular biology (\cite{Ful1, Ful2, Wh-Ba}). 

When the knot $K$ is given by $K=f(S^1)$ the {\em writhe} is given by the Gauss integral: 
\begin{equation}\label{def_writhe}
\textsl{Wr}\,(f)=\frac{\,1\,}{4\pi}\iint_{S^1\times S^1}\frac
{\det(f^{\p}(s),f^{\p}(t),f(s)-f(t))}{|f(s)-f(t)|^3}\,dsdt
\end{equation}
(the reader is referred to \cite{Al-Kl-Ta} for the details concerning writhe). 
Banchoff and White showed that the absolute value of the writhe is conformally invariant (\cite{Ban-Wh}). 
To be precise, they showed 
\begin{theorem}\label{thm_I_writhe}{\rm (\cite{Ban-Wh})} 
Suppose $I$ is an inversion in a sphere. 
Then we have 
\begin{eqnarray}\label{I_writhe}
\textsl{Wr}\,(I(K))=-\textsl{Wr}\,(K).
\end{eqnarray}
\end{theorem}
It is a corollary of 
\begin{theorem}\label{thm_I_twist}{\rm (\cite{Ban-Wh})} 
Suppose $I$ is an inversion in a sphere. 
Then we have 
\begin{eqnarray}\label{I_twist}
\textsl{Tw}(I(K))\equiv -\textsl{Tw}(K)\hspace{0.5cm}(\textrm{\rm mod}\>\> \mathbb Z). 
\end{eqnarray}
\end{theorem}

In this paper, using techniques in conformal geometry, we give new proofs of Theorem \ref{thm_I_writhe} and a special case of Theorem \ref{thm_I_twist} when $K$ has nowhere vanishing curvature and $\vect e_2$ is the unit principal normal vector field. 

\section{Notations}

Throughout the paper we use the following notations. 

We assume that a knot $K=f(S^1)$ is oriented. 
We denote the positive tangent vector $\dot f$ by $v=v(x)$. 

We denote the circle through $x,y$, and $z$ by $\Ga(x,y,z)$. 
When one of $x,y$, and $z$, say $z$ is $\infty$, $\Ga(x,y,\infty)$ means the line through $x$ and $y$. 
When $x$ is a point on a knot $K$, $\Ga(x,x,y)$ denotes the circle (or line) which is tangent to $K$ at $x$ that passes through $y$. 
We assume that it is oriented by $v$ at $x$. 
Especially, when $y=x$ $\Ga(x,x,x)$ denotes the oriented osculating circle. 
(We consider lines as circles through $\infty$.) 

The union of osculating circles $\displaystyle \bigcup_{x\in K}\Ga(x,x,x)$ is called the {\em curvature tube} (\cite{Ban-Wh}). 

Suppose $x$ is a point on a knot $K$. 
We denote a sphere through the osculating circle $\Ga(x,x,x)$ and $y$ $(y\ne x)$ by $\Si(x,x,x,y)$. 
It is uniquely determined generically, i.e. unless $y\in\Ga(x,x,x)$. 
When $y=\infty$ $\Si(x,x,x,\infty)$ means the plane through $\Ga(x,x,x)$. 

\section{Proof of Theorem \ref{thm_I_twist} in special case}
In this section we assume that a knot $K=f(S^1)$ has nowhere vanishing curvature and the unit normal vector field $\vect e_2$ along $K$ is given by the principal normal vectors. 
In this case the total twist is equal to $\frac1{2\pi}$ times the {\em total torsion} which we denote by $T\omega$: 
\[T\omega=\frac1{2\pi}\int_{K}\tau \,dx,\]
where $\tau$ denotes the torsion of the knot. 
We will show in this section 
\begin{proposition}\label{prop_thm2}
If both $K$ and $I(K)$ have nowhere vanishing curvatures then 
\begin{equation}\label{I_torsion}
T\omega(I(K))\equiv-T\omega(K) \hspace{0.5cm}(\textrm{\rm mod}\>\> \mathbb Z), 
\end{equation}
where $I$ is an inversion in a sphere. 
\end{proposition}

We remark that the above proposition can be proved by Theorem 6.3 of \cite{CSW} which gives 
\[T\omega (K)=\frac1{2\pi}\int_K\tau dx \equiv\frac1{2\pi}\int_KTd\rho \hspace{0.5cm}(\textrm{\rm mod}\>\> \mathbb Z),\]
where $T$ is the conformal torsion and $\rho$ is the conformal arc-length. 

\begin{lemma}\label{lem_1} 
Suppose a point $P$ does not belong to the osculating circle \newline
$\Ga(x_0, x_0, x_0)$ of a knot $K$ at $x_0$. 
{\rm (}We allow $P=\infty$.{\rm )} 
Then, infinitesimally speaking, the sphere $\Si(x, x, x, P)$ rotates around the circle $\Ga(x_0, x_0, P)$ at $x=x_0$ as $x$ travels in $K$. 
In other words, 
\begin{equation}\label{f_lem_1}
\displaystyle \lim_{x_1\to x_0}\left(\Si(x_0, x_0, x_0, P)\cap \Si(x_1, x_1, x_1, P)\right)\supset\Ga(x_0, x_0, P).
\end{equation}
\end{lemma}

\begin{proof} 
Let $\Si(x,y,z,w)$ denote a sphere through $x,y,z$, and $w$. 
Then 
\[\begin{array}{l}
\displaystyle \lim_{x_1\to x_0}\left(\Si(x_0, x_0, x_0, P)\cap \Si(x_1, x_1, x_1, P)\right)\\[2mm]
\displaystyle \>=\lim_{x_1\to x_0}\left(\lim_{y_0\to x_0, y_1\to x_1}\big(\Si(y_0, x_0, x_1, P)\cap \Si(x_0, x_1, y_1, P)\big)\right)\\[3mm]
\displaystyle \>\supset \lim_{x_1\to x_0}\Ga(x_0, x_1, P)=\Ga(x_0, x_0, P).
%
\end{array}\]

\medskip
There is an alternative computational proof using the Lorentzian exterior product introduced in \cite{La-OH}. 
\end{proof}

\medskip
We show that the total torsion of the image of $K$ by an inversion in a sphere with center $P$ is equal to the total angle variation of the rotation of the sphere $\Si(x, x, x, P)$ as $x$ goes around in $K$. 
In order to take into accound the sign of the torsion, we have to consider the orientations. 

\begin{figure}[htbp]
\begin{center}
\includegraphics[width=.8\linewidth]{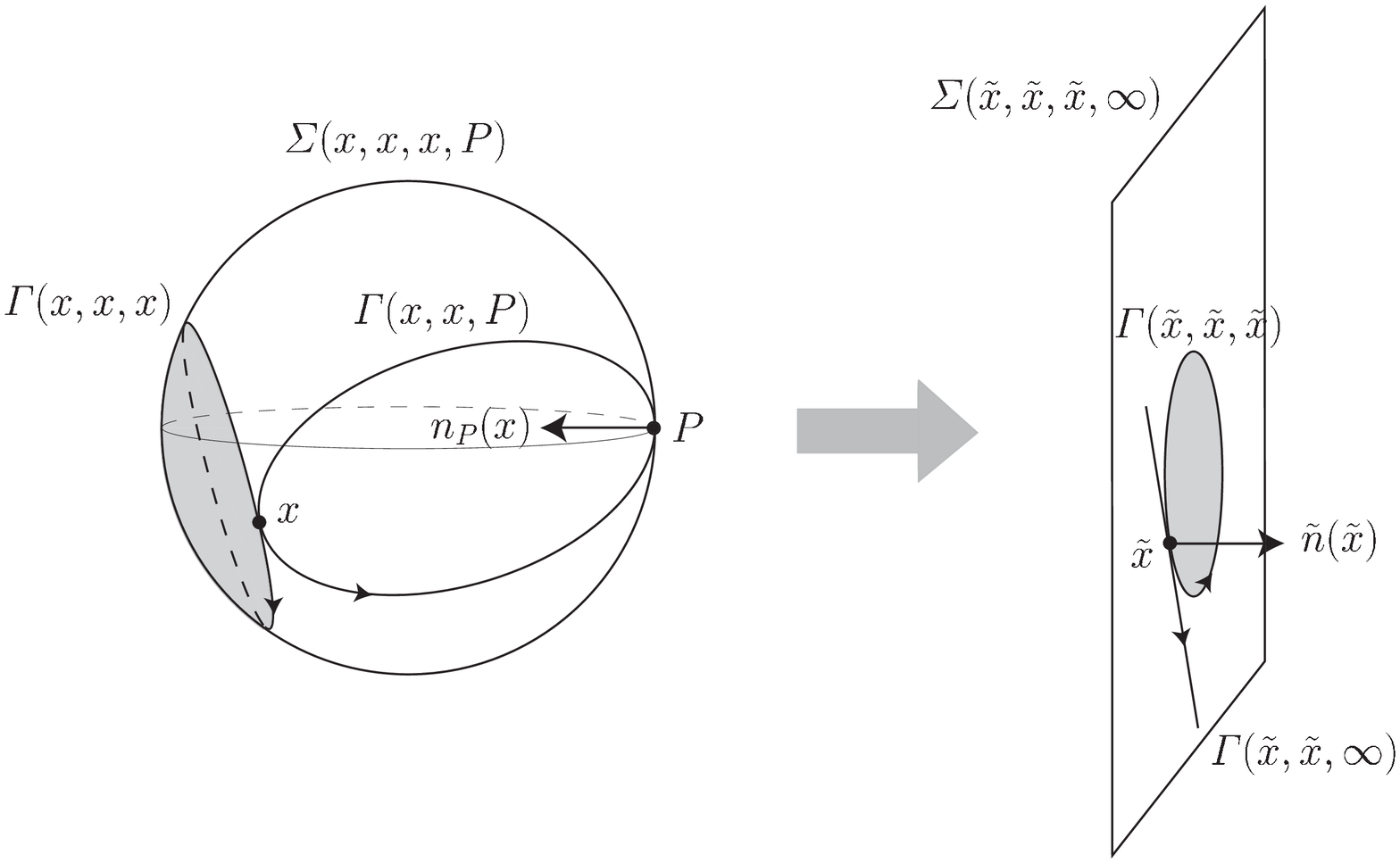}
\caption{}
\label{fig_tilde_n=-n_P}
\end{center}
\end{figure}

\smallskip
Fix $P\not\in\Ga(x,x,x)$, where we allow $P=\infty$. 
The osculating circle $\Ga(x,x,x)$ divide the sphere $\Si(x,x,x,P)$ into two domains. 
Let $D_1$ be one of the two that does not contain $P$ (gray disc of Figure \ref{fig_tilde_n=-n_P} left). 
Assume $\Si(x,x,x,P)$ is given the orientation such that the restriction to $D_1$ induces the same orientation to the boundary $\partial D_1$ as that of $\Ga(x,x,x)$ which is fixed in the previous section. 
Let $n(x)$ and $n_P(x)$ be the positive unit normal vectors to $\Si(x,x,x,P)$ at $x$ and $P$ respectively (Figure \ref{fig_tilde_n=-n_P} left). 
Let $\varPi(x)$ and $\varPi_P(x)$ be the normal planes to $\Ga(x,x,P)$ at $x$ and $P$ respectively. We assume that $\varPi(x)$ (or $\varPi_P(x)$) is oriented so that the algebraic intersection number of $\Ga(x,x,P)$ and $\varPi(x)$ (or $\varPi_P(x)$) at $x$ (or respectively, at $P$) is equal to $+1$. 

Since $\varPi(x)\perp\Ga(x,x,P)$
, Lemma \ref{lem_1} implies 
\begin{corollary}\label{cor_lemma_1}
Infinitesimally, the normal vector $n(x)$ to the sphere $\Si(x,x,x,P)$ rotates in the plane $\varPi(x)$ to $K$: 
\begin{equation}\label{n_varies_in_Pi}
\frac{d}{dx}n(x)\in \varPi(x).
\end{equation}
\end{corollary}
Let $I$ be an inversion in a sphere with center $P$. 
We denote $I(x)$ and $I(K)$ by $\tilde x$ and $\widetilde K$. 
Then $I$ maps the osculating circle $\Ga(x,x,x)$ to the osculating circle $\Ga(\tilde x, \tilde x, \tilde x)$ of $\widetilde K$ at $\tilde x$, and the sphere $\Si(x,x,x,P)$ to the plane $\Si(\tilde x, \tilde x, \tilde x, \infty)$. 
Let $\tilde n(\tilde x)$ be the positive unit normal vector to $\Si(\tilde x, \tilde x, \tilde x, \infty)$. 
Then 
\begin{equation}\label{f_tilde_n=-n_P}
\tilde n(\tilde x)=-n_P(x)
\end{equation}
(Figure \ref{fig_tilde_n=-n_P}). 

The above convention of orientation implies that $\tilde n(\tilde x)$ is equal to the unit binormal vector of $\widetilde K$. 
Lemma \ref{lem_1} implies that $\tilde n(\tilde x)$ rotates in $\varPi({\tilde x})$ at $\tilde x$. 
Our convention of the orientation of $\varPi({\tilde x})$ implies 
\begin{lemma}\label{lem_2}
The torsion of $\widetilde K$ is equal to the angle velocity of $\tilde n$ with respect to the arc-length $\tilde s$: 
\begin{equation}\label{f_torsion_inverted_knot}
\tau(\widetilde K)(\tilde x)=\e(\tilde x)\,\left|\wfrac{d\tilde n}{d\tilde s}(\tilde x)\right| \hspace{0.5cm}(\e(\tilde x)\in\{+1, -1\}),
\end{equation}
where $\e(\tilde x)$ is the signature of the rotation of $\tilde n$ at $\tilde x$ with respect to the orientation of $\varPi({\tilde x})$. 
\end{lemma}
In other words, we have 
\[\tau(\widetilde K)(\tilde x)=\tilde v\cdot \left( \tilde n \times \frac{d\tilde n}{d\tilde s}\right)\,, 
\]
where $\tilde v$ 
denotes the positive unit tangent vector to $\widetilde K$. 
\begin{figure}[htbp]
\begin{center}
\includegraphics[width=.5\linewidth]{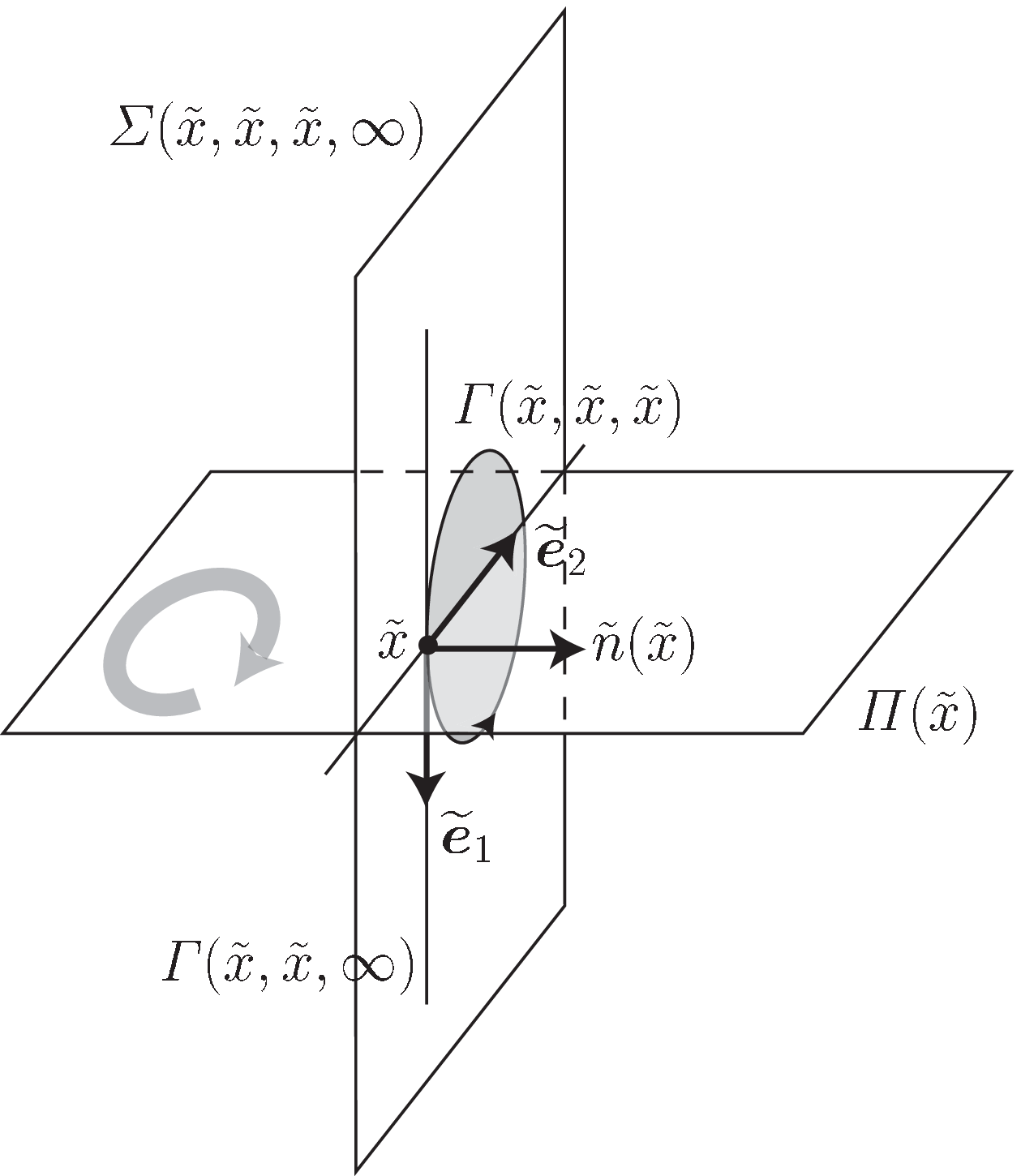}
\caption{}
\label{fig_lem_2}
\end{center}
\end{figure}

\begin{figure}[htbp]
\begin{center}
\includegraphics[width=.6\linewidth]{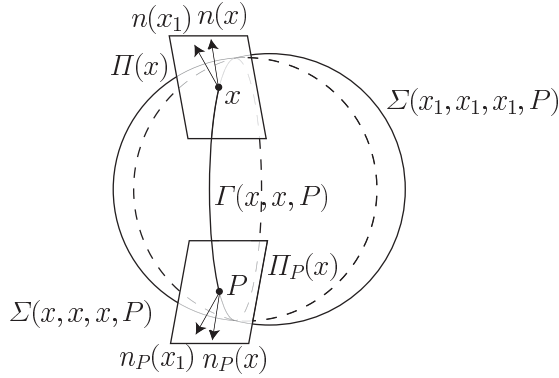}
\caption{$\angle n(x)\cdot n(x_1)=\angle n_P(x)\cdot n_P(x_1)=\angle \Si(x,x,x,P)\cdot\Si(x_1, x_1, x_1, P)$}
\label{fig_n_varies_in_Pi}
\end{center}
\end{figure}

Since the positive unit tangent vector to $\Ga(x,x,P)$ at $P$ is equal to $-\tilde v(\tilde x)$, the orientations of $\varPi_P(x)$ is opposite to that of $\varPi({\tilde x})$. 
Therefore, (\ref{f_tilde_n=-n_P}) and Lemma \ref{lem_2} imply 
\[\tau(\widetilde K)(\tilde x)=-\e_P\left|\wfrac{ds}{d\tilde s}\right|\,\left|\wfrac{dn_P}{ds}\right| \hspace{0.5cm}(\e_P\in\{+1, -1\}),
\]
where $\e_P$ is the signature of the rotation of $n_P$ with respect to the orientation of $\varPi_P({x})$. 
Then Lemma \ref{lem_1} implies $\displaystyle \left|\wfrac{dn_P}{ds}\right|=\left|\wfrac{dn}{ds}\right|$ (Figure \ref{fig_n_varies_in_Pi}), therefore we have 
\begin{lemma}\label{lem_3}
The torsion of $\widetilde K$ at $\tilde x$ is equal to the negative of the angle velocity of $n$ in $\varPi(x)$ with respect to the arc-length $s$ of $K$ 
up to the multiplication by the Jacobian: 
\begin{equation}\label{f_lem_3}
\tau(\widetilde K)(\tilde x)=-\e(x)\left|\wfrac{ds}{d\tilde s}\right|\left|\wfrac{dn}{ds}(x)\right| \hspace{0.5cm}(\e(x)\in\{+1,-1\}),
\end{equation}
where $\e(x)$ is the signature of the rotation of $n$ at $x$ with respect to the orientation of $\varPi({x})$. 
\end{lemma}
In other words, we have 
\[\tau(\widetilde K)(\tilde x)=-\left|\wfrac{ds}{d\tilde s}\right|\,v\cdot \left(n\times\frac{dn}{ds}\right)\,.
\]
Thus we are led to 
\begin{proposition}\label{prop_4}
Let $K$ be a knot with nowhere vanishing curvature, $I$ an inversion in a sphere with center $P$, and $\widetilde K=I(K)$. 
Assume that $\widetilde K$ has nowhere vanishing curvature, which happens if and only if $P$ is not contained in the curvature tube of $K$. 
Then the total torsion of $\widetilde K$ is the negative of the total angle variation of the positive unit normal vector $n(x)$ to the sphere $\Si(x, x, x, P)$ as $x$ goes around in $K$: 
\begin{equation}\label{f_prop_4}
\int_{\widetilde K}\tau(\widetilde K)(\tilde x)d\tilde x
=-\int_{K}\e\left|{dn}\right|
=-\int_{K}v\cdot (n\times{dn})\,. 
\end{equation}
\end{proposition}
\begin{corollary}\label{cor_5}
Let $K$ be a knot with nowhere vanishing curvature, and $I_j$ $(j=1,2)$ an inversion in a sphere with center $P_j$ which is not contained in the curvature tube of $K$. 
Then the total torsion of $I_1(K)$ and $I_2(K)$ coincide modulo $2\pi\mathbb Z$: 
\[\int_{I_1(K)}\tau d\tilde x_1\equiv \int_{I_2(K)}\tau d\tilde x_2 \hspace{0.5cm}(\textrm{\rm mod}\>\> 2\pi\mathbb Z).\]
In other words, 
\[T\omega(I_1(K))-T\omega(I_2(K))\in\mathbb Z.\]
\end{corollary}
\begin{proof}
Let $n_j$ $(j=1,2)$ be the positive unit normal vector to $\Si(x,x,x,P_j)$ at $x$, and $\theta_{21}$ the angle from $n_1$ to $n_2$ in the oriented plane $\varPi(x)$. 
Then Corollary \ref{cor_lemma_1} and the above Proposition imply 
\[\int_{I_1(K)}\tau d\tilde x_1-\int_{I_2(K)}\tau d\tilde x_2=\int_Kd\theta_{21},\]
which is equal to $2\pi k$ $(k\in\mathbb Z)$ since $K$ is closed. 
\end{proof}

\medskip
\begin{proofofproposition} {\bf \ref{prop_thm2}}: 
Put $P_2=\infty$ in the above Corollary. 
\end{proofofproposition}

\section{Proof of Theorem \ref{thm_I_writhe}}
We first prove the following fact in our context, which implies Theorem \ref{prop_thm2} in the case when the knot has nowhere vanishing curvature. 
\begin{proposition}\label{prop_thm1} {\rm (\cite{Ca1}, \cite{Po1})} 
Suppose $K$ has nowhere vanishing curvature. 
Then 
\[4\pi \textsl{Wr}\,(K)=-2\int_K\tau +4\pi k \hspace{0.4cm}\mbox{for some}\>\>k\in\mathbb Z.\]
In other words, $\textsl{Wr}+T\omega$ is an integer. 
\end{proposition}
\begin{proof}
Assume $K=f(S^1)$ is parametrized by the arc-length. 
Suppose $S^1=[0, L]/\sim$, where $L$ is the length of $K$. 
Let $\varphi$ be a map from $S^1\times S^1\setminus\Delta$ to $S^2$ given by 
\[\varphi(s,t)=\frac{f(s)-f(t)}{|f(s)-f(t)|}\hspace{0.5cm}(s\ne t).\]
Then the integrand of the writhe is equal to the pull-back of the standard area element of $S^2$ by $\varphi$: 
\[\frac{\det(f^{\p}(s),f^{\p}(t),f(s)-f(t))}{|f(s)-f(t)|^3}=\varphi\cdot(\varphi_s\times\varphi_t),\]
where $\varphi_s=\nekofrac{\partial \varphi}{\partial s}$ and $\varphi_t=\nekofrac{\partial \varphi}{\partial t}$. 
Therefore, $4\pi$ times the writhe is equal to the signed area of the image of 
\[D=\{(s,t)\,|\,0\le s\le L, s<t<s+L\}\]
by $\varphi$. 
Let $S$ be the closure of $\varphi(D)$. 
Then $S$ is an oriented ``{\sl surface}" (a continuous image of $S^1\times [0,L]$) possibly with self-overlaps in $S^2$ whose ``{\sl boundary}" (the image of $S^1\times \{0,L\}$) is given by 
\[\partial S=C_+\cup (-C_-),\]
where $C_+$ and $C_-$ denote the positive and negative tangential indicatrices: \setlength\arraycolsep{1pt}
\[\begin{array}{rl}
C_+=&\varphi(\Delta_+)=\{v(s)=\dot f(s)\,|\,0\le s\le L\},\\[2mm]
C_-=&\varphi(\Delta_-)=\{-v(s)=-\dot f(s)\,|\,0\le s\le L\}. 
\end{array}\]
We assume that the indicatirices are oriented so that $s$ increases in the positive direction. 

Let $\Ga(s)$ denote a great circle in $S^2$ which is tangent to $C_+$ at $v(s)$ (and to $C_-$ at $-v(s)$). 
We remark that $C_+$ has nowhere vanishing tangent vector since the knot has nowhere vanishing curvature. 
Note that $C_+\cup C_-$ is the envelope of the family $\{\Ga(s)\}_{0\le s\le L}$. 
We assume that $\Ga(s)$ has an orientation compatible with that of $C_+$ at $v(s)$. 
Let $\Ga_+(s)$ be a semi-circle of $\Ga(s)$ from $v(s)\in C_+$ to $-v(s)\in C_-$ in the positive direction of $\Ga(s)$, and $S^{\p}$ a region in $S^2$ swept by $\Ga_+(s)$ as $s$ varies in $[0,L]$. 
Then $S^{\p}$ is given by 
\[S^{\p}=\{w(s,t)=(\cos t)v(s)+(\sin t)\frac{\dot v}{|\dot v|}(s)\,|\,0\le s\le L, \,0\le t\le \pi\}.\]
It is an oriented ``{\sl surface}" (a continuous image of $S^1\times [0,\pi]$) possibly with self-overlaps in $S^2$ whose ``{\sl boundary}" (the image of $S^1\times \{0,\pi\}$) is $C_+\cup (-C_-)$. 
The signed area of $S^{\p}$ is given by 
\[\textsl{Area}\,(S^{\p})=\int_0^L\int_0^{\pi}w\cdot(w_s\times w_t)dsdt. \]
Since $\dot v/|\dot v|$ is equal to the principal normal vector $\vect e_2$ of the knot $K$, we have 
\[\begin{array}{rcl}
w\cdot(w_s\times w_t)
&=&\det\left((\cos t)\vect e_1+(\sin t)\vect e_2, (\cos t)\kappa \vect e_2+(\sin t)\dot{\vect e_2}, -(\sin t)\vect e_1+(\cos t)\vect e_2\right)\\[1mm]
&=&-\tau \sin t,\end{array}\]
which implies 
\[\textsl{Area}\,(S^{\p})=-2\int_0^L\tau \,ds=-4\pi \,T\omega(K).\]

Since both $S$ and $S^{\p}$ have the boundary $C_+\cup (-C_-)$, $S\cup (-S^{\p})$ is a cycle of $S^2$, i.e. an oriented ``{\sl surface}" (a continuous image of a torus) possibly with self-overlaps without a boundary. 
Therefore, the signed area of $S\cup (-S^{\p})$ is equal to $4\pi k$ for some integer $k$. 
It follows that 
\[\textsl{Area}\,(S\cup (-S^{\p}))=\textsl{Area}\,(S)-\textsl{Area}\,(S^{\p})=4\pi \textsl{Wr}+4\pi T\omega=4\pi k \>\>(k\in\mathbb Z), \]
which completes the proof. 
\end{proof}

\begin{corollary}\label{cor_prop_thm1} 
If both $K$ and $I(K)$ have nowhere vanishing curvatures then 
\[\textsl{Wr}\,(I(K))=-\textsl{Wr}\,(K).\]
\end{corollary}

\begin{proof} 
Let $I_j$ $(j=0,1)$ be an inversion in a sphere with center $P_j$ which is not contained in the curvature tube of $K$. 
Then Proposition \ref{prop_thm1} implies that 
\[\begin{array}{r}
\textsl{Wr}\,(I_0(K))+T\omega(I_0(K))\in\mathbb Z, \\[1mm]
\textsl{Wr}\,(I_1(K))+T\omega(I_1(K))\in\mathbb Z. 
\end{array}\]
Then Corollary \ref{cor_5} implies 
\[\textsl{Wr}\,(I_0(K))-\textsl{Wr}\,(I_1(K))\in\mathbb Z.\]
Join $P_0$ and $P_1$ by a smooth path $P_t$. 
Let $I_t$ be an inversion in a sphere with center $P_t$. 
Then $\textsl{Wr}\,(I_t(K))$ is a continuous function of $t$ (\cite{Po1}), and hence 
\[\textsl{Wr}\,(I_0(K))=\textsl{Wr}\,(I_1(K)).\]
When $P_1$ goes to $\infty$ and the radius of the sphere of the inversion also goes to $+\infty$, $I_1(K)$ approaches the mirror image of $K$ and hence $\textsl{Wr}\,(I_1(K))$ approaches $-\textsl{Wr}\,(K)$, which completes the proof. 
\end{proof}

\medskip
\begin{proofoftheorem} {\bf \ref{thm_I_writhe}}: 
Suppose the curvature of $K$ vanishes somewhere. 
We have only to show that $K$ can be approximated, with respect to the $C^2$-topology, by a knot with non-vanishing curvature. 
This can be done as follows. 
The curvature tube of $K$ is non-compact as it contains a line. 
But we can still find a point $P$ with a very big distance from $K$ which is not contained in the cuvature tube. 
Let $I$ be an inversion in a sphere with center $P$ and radius approximately equal to the distance between $P$ and $K$. 
We can get a desired knot by taking the mirror image of $I(K)$ thus constructed. 
\end{proofoftheorem}

\end{document}